\documentclass[12pt]{amsart}
\usepackage{amscd}
\usepackage{amsmath}
\usepackage{amssymb}
\usepackage{units}
\usepackage[all]{xy}
\def\frk{\frak}               

\def\Phi{{\frk n}}
\def\Phi{{\frk N}}
\def\opn#1#2{\def#1{\operatorname{#2}}} 
\opn\chara{char} \opn\length{\ell} \opn\pd{pd} \opn\rk{rk}
\opn\projdim{proj\,dim} \opn\injdim{inj\,dim} \opn\rank{rank}
\opn\depth{depth} \opn\sdepth{sdepth} \opn\fdepth{fdepth}
\opn\grade{grade} \opn\height{height} \opn\embdim{emb\,dim}
\opn\codim{codim}  \opn\min{min} \opn\max{max}

\opn\Tr{Tr} \opn\bigrank{big\,rank}
\opn\superheight{superheight}\opn\lcm{lcm}
\opn\trdeg{tr\,deg}
\opn\reg{reg} \opn\lreg{lreg} \opn\ini{in} \opn\lpd{lpd}
\opn\size{size}
\opn\div{div} \opn\Div{Div} \opn\cl{cl} \opn\Cl{Cl}
\opn\Spec{Spec} \opn\Supp{Supp} \opn\supp{supp} \opn\Sing{Sing}
\opn\Ass{Ass} \opn\Min{Min}
\opn\Ann{Ann} \opn\Rad{Rad} \opn\Soc{Soc}
\opn\Im{Im} \opn\Ker{Ker} \opn\Coker{Coker} \opn\Am{Am}
\opn\Hom{Hom} \opn\Tor{Tor} \opn\Ext{Ext} \opn\End{End}
\opn\Aut{Aut} \opn\id{id}  \opn\deg{deg}

\opn\nat{nat}
\opn\pff{pf}
\opn\Pf{Pf} \opn\GL{GL} \opn\SL{SL} \opn\mod{mod} \opn\ord{ord}
\opn\Gin{Gin} \opn\Hilb{Hilb}
\opn\aff{aff} \opn\con{conv} \opn\relint{relint} \opn\st{st}
\opn\lk{lk} \opn\cn{cn} \opn\core{core} \opn\vol{vol}
\opn\link{link} \opn\star{star}
\opn\gr{gr}

\def\pot#1#2{#1[\kern-0.28ex[#2]\kern-0.28ex]}

\opn\dirlim{\underrightarrow{\lim}}
\opn\inivlim{\underleftarrow{\lim}}
\let\Dirsum=\bigoplus

\let\to=\rightarrow

\def\Implies{\ifmmode\Longrightarrow \else
        \unskip${}\Longrightarrow{}$\ignorespaces\fi}
\def\implies{\ifmmode\Rightarrow \else
        \unskip${}\Rightarrow{}$\ignorespaces\fi}
\def\iff{\ifmmode\Longleftrightarrow \else
        \unskip${}\Longleftrightarrow{}$\ignorespaces\fi}

\let\:=\colon
\newtheorem{Theorem}{Theorem}[section]
\newtheorem{Lemma}[Theorem]{Lemma}
\newtheorem{Corollary}[Theorem]{Corollary}
\newtheorem{Proposition}[Theorem]{Proposition}
\newtheorem{Remark}[Theorem]{Remark}

\newtheorem{Example}[Theorem]{Example}

\newtheorem{Definition}[Theorem]{Definition}

\let\epsilon\varepsilon
\let\phi=\varphi
\let\kappa=\varkappa
\def\qed{\ifhmode\textqed\fi
      \ifmmode\ifinner\quad\qedsymbol\else\dispqed\fi\fi}
\def\textqed{\unskip\nobreak\penalty50
       \hskip2em\hbox{}\nobreak\hfil\qedsymbol
       \parfillskip=0pt \finalhyphendemerits=0}
\def\dispqed{\rlap{\qquad\qedsymbol}}

\opn\dis{dis}
\def\pnt{{\raise0.5mm\hbox{\large\bf.}}}

\opn\Lex{Lex}

\begin{document}

\title{\bf Graph and depth of a monomial squarefree ideal}

\author{ Dorin Popescu }

\thanks{The  Support from  the CNCSIS grant PN II-542/2009 of Romanian Ministry of Education, Research and Inovation is gratefully acknowledged.}

\address{Dorin Popescu, Institute of Mathematics "Simion Stoilow",
University of Bucharest, P.O.Box 1-764, Bucharest 014700, Romania}
\email{dorin.popescu@imar.ro}

\maketitle
\begin{abstract} Let $I$ be a monomial squarefree ideal of a polynomial ring $S$ over a field $K$ such that the sum of every three different of its minimal prime ideals is the maximal ideal of $S$, or more general a constant ideal. We associate to $I$ a graph on $[s]$, $s=|\Min S/I|$ on which we may {\em read} the depth of $I$. In particular,
 $\depth_S\ I$ does not depend of char $K$. Also we show that $I$ satisfies the Stanley's Conjecture.

  \vskip 0.4 true cm
 \noindent
\end{abstract}

\section*{Introduction}

 Let  $S=K[x_1,\ldots,x_n]$, $n\in {\bf N}$ be a polynomial ring over a field $K$, and  $I\subset S$  a monomial squarefree ideal  with  minimal prime ideals  $P_1,\ldots,P_s$ (here we study only the  monomial squarefree ideals). After \cite{L} the size of $I$ is the number $v+(n-h)-1$, where $h$ is the height of $\sum_{j=1}^sP_j$ and $v$ is the minimal number $e$ for which there exist integers $i_1<i_2<\cdots <i_e$ such that $\sum_{k=1}^eP_{i_k}=\sum_{j=1}^sP_j$. Similarly, we defined in \cite{P2} the bigsize of $I$, which is the number  $t+(n-h)-1$, where $t$ is the minimal number $e$ such that for all integers $i_1<i_2<\cdots <i_e$ it holds $\sum_{k=1}^eP_{i_k}=\sum_{j=1}^sP_j$. Clearly $bigsize(I)\geq size(I)$.
 Lyubeznik \cite{L} showed that $\depth I\geq  1+\size I$.

  Let $I\subset S$ be a monomial ideal of $S$, $u\in I$ a monomial and $uK[Z]$, $Z\subset \{x_1,\ldots ,x_n\}$ the linear $K$-subspace of $I$ of all elements $uf$, $f\in K[Z]$. A  presentation of $I$ as a finite direct sum of such spaces ${\mathcal D}:\
I=\Dirsum_{i=1}^ru_iK[Z_i]$ is called a Stanley decomposition of $I$. Set $\sdepth
(\mathcal{D})=\min\{|Z_i|:i=1,\ldots,r\}$ and
\[
\sdepth\ I :=\max\{\sdepth \ ({\mathcal D}):\; {\mathcal D}\; \text{is a
Stanley decomposition of}\;  I \}.
\]

The Stanley's Conjecture \cite{S} says that $\sdepth\ I\geq \depth\ I$.  This conjecture holds for arbitrary  monomial squarefree ideals if $n\leq 5$
by \cite{P} (see especially the arXiv version), or for intersections of four monomial
prime ideals by \cite{AP}, \cite{P2}. In the case of non squarefree monomial ideals $J$ an important inequality is $\sdepth J\leq \sdepth \sqrt{J}$ (see \cite[Theorem 2.1]{Is}).
Similarly to Lyubeznik's result, it holds  $\sdepth I\geq  1+\size I$ by \cite[Theorem 3.1]{HPV}. If $bigsize(I)= size(I)$
then $I$ satisfies the Stanley's Conjecture by \cite[Theorems 1.2, 3.1]{HPV}.

The purpose of this paper is to study the case when $bigsize(I)=2$, $ size(I)=1$. In the case $\sum_{j=1}^sP_j=m=(x_1,\ldots,x_n)$, we associate to $I$ a graph $\Gamma$ on $[s]$ given by  {\em $\{ij\}$ is an edge  if and only if $P_i+P_j=m$}. We express the depth of $I$ in terms of the properties of $\Gamma$ and  of
$q(I)=\min\{\dim S/(P_i+P_j): j\not =i,P_i+P_j\not =m\}.$ We note that \cite[Lemmas 3.2, 3.2]{P2} say, in particular, that $\depth_S\ I=2$ if and only if  $\Gamma$
is a join graph. Our Corollary \ref{carac} says that  if $q(I)>1$ then $\depth_S\ I=2+q(I)$  if and only if  $\Gamma$ is a so called concatenation of several graphs on two vertices having no edges. Thus knowing $q(I)$, $\depth_S\ I$ can be {\em read} on $\Gamma$(see Corollary \ref{car}). It follows that for  a monomial squarefree ideal $I\subset S$   such that the sum of every three different of its minimal prime ideals is  a constant ideal (for example $m$),   $\depth_S\ I$ does not depend of char $K$ (see Theorem \ref{char}) and the Stanley's Conjecture holds (see Theorem \ref{m}).

It is well known that $\depth_S\ I$ depends of the characteristic of $K$ if $bigsize(I)=3$, $size(I)=2$ (see our Remark \ref{f}), so it is very likely that this case is much harder for proving Stanley's Conjecture. Several people ask if there exist  examples when the special Stanley decomposition of \cite{AP}, \cite{P2}, or the splitting variables in the terminology of \cite{HPV} do not help in proving Stanley's Conjecture since there exists no {\em good} main prime ideal. Our Example \ref{ex3} is such an example.

\vskip 1 cm
\section{Depth  two and three}

Let $S=K[x_1,\ldots,x_n]$, $n\in {\bf N}$ be a polynomial ring over a field $K$ and ${\tilde S}=K[x_1,\ldots,x_{n-1}]\subset S$.
We start  reminding the following two lemmas from \cite{P}.

\begin{Lemma}\label{good}
Let  $I, J\subset {\tilde S}$, $I\subset J$, $I\not=J$ be two monomial
ideals, $T=(I+x_nJ)S$ such that
\begin{enumerate}
\item{} $\depth_{{\tilde S}}\ {\tilde S}/I\not =\depth_S\ S/T -1,$
\item{} $\sdepth_{\tilde S}\ I\geq \depth_{\tilde S}\ I,$ \ \ \
$\sdepth_{\tilde S}\ J\geq \depth_{\tilde S}\ J.$
\end{enumerate}
Then $\sdepth_{S}\ T\geq \depth_{S}\ T.$
\end{Lemma}

\begin{Lemma}\label{bad} Let $I, J\subset {\tilde S}$, $I\subset J$,
$I\not=J$ be two monomial ideals, $T=(I+x_nJ)S$ such that
\begin{enumerate}
\item{} $\depth_{{\tilde S}}\ {\tilde S}/I=\depth_S\ S/T -1,$
\item{} $\sdepth_{{\tilde S}}\ I\geq \depth_{{\tilde S}}\ I,$
\item{} $\sdepth_{{\tilde S}}\ J/I\geq \depth_{{\tilde S}}\ J/I.$
\end{enumerate}
Then $\sdepth_{S}\ T\geq \depth_{S}\ T.$
\end{Lemma}
The above lemmas allow us to show Stanley's Conjecture in a special case.
\begin{Proposition}\label{dep2}
Let $T\subset S$ be a monomial squarefree ideal. If $S/T$ is Cohen-Macaulay of dimension $2$ then  $\sdepth_{S}\ T\geq \depth_{S}\ T.$
\end{Proposition}
\begin{proof}
We use induction on $n$, case $n\leq 5$ being given in \cite{P}. Suppose $n>5$. Then $T$ has the form $T=I+x_nJ$ for two monomial squarefree ideals
$I,J\subset {\tilde S}$, in fact $I=T\cap {\tilde S}$, $J=(T:x_n)\cap {\tilde S}$.
 Note that $\dim {\tilde S}/I=\dim S/(T,x_n)\leq 2$ and $\dim S/JS=\dim ((x_n)+T)/T\leq 2$ and so $\depth_{{\tilde S}}\ {\tilde S}/I\leq 2$,
 $\depth_{{\tilde S}}\ {\tilde S}/J\leq 1$. If $\depth_{{\tilde S}}\ {\tilde S}/I= 2$ then $\sdepth_{{\tilde S}}\ I\geq\depth_{{\tilde S}}\ I$
 by induction hypothesis. If $\depth_{{\tilde S}}\ {\tilde S}/I= 1$ (by \cite[Proposition 1.2]{R} $\depth_{{\tilde S}}\ {\tilde S}/I>0$) then $\depth_{{\tilde S}}\ I=2=1+size(I)\leq
 \sdepth_{{\tilde S}}\ I$ by \cite[Theorem 3.1]{HPV} and similarly for $J$.
 As $\dim\ J/I\leq \dim\ {\tilde S}/I\leq \dim\ S/T=2$ we have  $\sdepth_{{\tilde S}}\ J/I\geq \depth_{{\tilde S}}\ J/I$ by \cite{Po}. Now the result is a consequence of the Lemmas \ref{good}, \ref{bad} if $I\not=J$, otherwise $T=IS$ and we may apply \cite[Lemma 3.6]{HPV}.
\end{proof}

 Let $I=\cap_{i=1}^s P_i$, $s\geq 3$  be the intersection of the minimal monomial prime ideals of $S/I$.    Assume that $\Sigma_{i=1}^s P_i=m$
and the bigsize of $I$ is two. Set $$q=q(I)=\min\{\dim S/(P_i+P_j): j\not =i,P_i+P_j\not =m\}.$$

 We will need the following two lemmas from \cite{P2}.

\begin{Lemma} \label{g2} If $P_1+P_2\not =m$ and $P_k+P_e=m$ for all  $k,e>2$, $k\not =e$ then
\begin{enumerate}
\item{} $\depth_S S/I\in \{1,2,1+q\}$,
\item{} $\depth_S S/I=1$ if and only if there exists $j>2$ such that $P_1+P_j=m=P_2+P_j$,
\item{}  $\depth_S S/I>2$ if and only if $q>1$ and each $j>2$ satisfies either $$P_1+P_j\not=m =P_2+P_j,\ \mbox{or}$$
$$P_2+P_j\not=m =P_1+P_j,$$
\item{}  $\depth_S S/I=2$ if and only if  the following conditions hold:
\begin{enumerate}
\item{}  each $j>2$ satisfies either $P_1+P_j\not=m$ or
$P_2+P_j\not=m,$
\item{} $q=1$  or there exists a $k>2$ such that  $$P_1+P_k\not=m \not=P_2+P_k,$$
\end{enumerate}
\item{} $\sdepth_S I\geq \depth_S I$.

\end{enumerate}
\end{Lemma}

\begin{Lemma}\label{b2} Suppose that whenever there exist $i\not =j$ in $[s]$ such that $P_i+P_j\not =m$  there exist also $k\not =e$ in $[s]\setminus\{i,j\}$
such that $P_k+P_e\not =m$
 (that is the complementary case of  the above lemma). Then
\begin{enumerate}
\item{} $\depth_S S/I\in \{1,2,1+q\}$.
\item{}  $\depth_S S/I=1$ if and only if  after a renumbering of $(P_i)$ there exists $1\leq c<s$ such that $P_i+P_j=m$ for each $c<j\leq s$ and $1\leq i\leq c$.
\end{enumerate}
\end{Lemma}
These two lemmas allow us to show the following useful proposition.

\begin{Proposition}\label{1} Suppose that $P_1=(x_1,\ldots,x_r)$, $1\leq r<n$, $S'=K[x_{r+1},\ldots,x_n]$ and $P_1+P_2\not =m\not =P_1+P_3$, $P_2+P_3\not =m$. Then $depth_S\ S/I\leq 2$, in particular
$$\sdepth_{S'} (P_2\cap P_3\cap S')\geq 2\geq \depth_S\ S/I.$$
\end{Proposition}
\begin{proof}
Apply induction on $s$, the cases $s=3,4$ follows from \cite{AP}, \cite{P2}. Suppose that $s>4$. Set $E=S/(P_1\cap  P_3\cap \ldots \cap P_s)\oplus S/(P_1 \cap P_2\cap P_4\cap \ldots\cap P_s)$  and\\
$F=S/(P_1\cap  (P_2+P_3)\cap P_4\cap\ldots \cap P_s)$. Note that if $P_i\subset P_2+P_3$ for some $i\not =2,3$ then $P_2+P_3=P_i+P_2+P_3=m$ because
bigsize of $I$ is two. Contradiction! Thus the bigsize of $F$ is one and so $\depth_S S/F=1$ by \cite{P2}. From the following exact sequence
$$0\to S/I\to E\to F\to 0$$
we get $\depth_S S/I=2$ if $\depth_SE>1$. Otherwise, suppose that $G=S/(P_1 \cap P_2\cap P_4\cap \ldots\cap P_s)$ has depth one. Then after renumbering $(P_i)$ we may suppose that there exists $c\not =3$, $1\leq c<s$ such that $P_i+P_j=m$ for all $1\leq i\leq c$, $c<j\leq s$, $i,j\not =3$ (see Lemmas \ref{g2}, \ref{b2}). In fact we may renumber only $(P_e)_{e>3}$ and take $c>3$ because $P_1+P_2\not =m$.
Set $M=S/P_1\cap\ldots\cap P_c$ and $N=M\oplus S/P_3\cap P_{c+1}\cap\ldots P_s$. In the following exact sequence
$$ 0\to S/I\to N\to S/P_3\to 0$$
we have  the depth of all modules $\leq \depth_S S/P_3$. By Depth Lemma \cite{Vi} it follows $\depth_S S/I=\depth_S N$ and so  $\depth_S S/I\leq\depth_S M$. Applying the induction hypothesis we get $\depth_S M\leq 2$, that is $\depth_S S/I\leq 2$.
Finally, by \cite{PQ} we have
  $$\sdepth_{S'} (P_2\cap P_3\cap S')\geq \depth_{S'} (P_2\cap P_3\cap S')=1+\depth_{S'} S'/(P_2\cap P_3\cap S')=$$
  $$1+\depth_S S/(P_1+P_2)\cap (P_1+P_3)=2$$ because $P_1+P_2+P_3=m$.
\end{proof}

\begin{Corollary}\label{ind} Suppose that $bigsize(I)=size(I)\leq 2$. Then $\depth_S\ I$ does not depend on the characteristic of $K$
\end{Corollary}
\begin{proof} If $bigsize(I)=size(I)=1$ then $\depth_S\ I=2$ by \cite[Corollary 1.6]{P2} and so does not depend on the characteristic of $K$. If $bigsize(I)=size(I)=2$ then $\depth_S\ I\leq 3$ by the above proposition and so $\depth_S\ I= 3$ by \cite{L} independently of char $K$.
\end{proof}

 \begin{Theorem} \label{dep3}
 If $\depth_S\ I\leq 3$ then $\sdepth_S\ I\geq  \depth_S\ I$.
 \end{Theorem}
 \begin{proof} By \cite{HPV} we have $\sdepth_S\ I\geq 1+size(I) \geq 2$ and it is enough to consider the case $\depth_S\ I=3$, that is  $\depth_S\ S/I=2$.  If
 $\dim\ S/I=2$ then we may apply  Proposition \ref{dep2}, otherwise we may suppose that $\dim\ S/P_i\geq 3$ for an $i$, let us say $i=1$. We may suppose that $P_1=(x_1,\ldots,x_r)$ for some $r<n$, thus $n\geq r+3$. Set $S''=K[x_1,\ldots,x_r]$, $S'=K[x_{r+1},\ldots,x_n]$.

 Applying \cite[Theorem 1.5]{P2} for ${\mathcal F}$ containing some $\tau_j=\{ j\}$,  and $\tau_{jk}=\{ j,k\}$  $1<j,k\leq s$,$j\not=k$
 we get $\sdepth_S I\geq \min\{A_0,\{A_{\tau_j}\}_{\tau_j\in {\mathcal F}}\},\{A_{\tau_{jk}}\}_{\tau_{jk}\in {\mathcal F}}\}$ for
$A_0=\sdepth_S(I\cap S'')S$ if $I\cap S''\not=0$ or $A_0=n$ otherwise, and
$$A_{\tau}\geq \sdepth_{S_{\tau}}J_{\tau} +\sdepth_{S'}L_{\tau},$$
where $J_{\tau}=\cap_{e\not\in\tau}P_e\cap S_{\tau} \not =0$,  $S_{\tau}=K[\{x_u:x_u\in S'',x_u\not \in \Sigma_{e\in \tau}P_e\}]$, $L_{\tau}=\cap_{e\in \tau}(P_e\cap S')\not =0$.
If $P_1+P_j\not =m$ then $$A_{\tau_j}\geq \sdepth_{S_{\tau_j}}J_{\tau_j} +\sdepth_{S'}(P_j\cap S')\geq 1+\dim S/(P_1+P_j)+\lceil\frac{\height (P_j\cap S')}{2}\rceil,$$
where $\lceil a\rceil$, $a\in {\bf Q}$ denotes the smallest integer $\geq a$.
Thus $A_{\tau_j}\geq 3= \depth_S I$. If $P_1+P_j =m$ then $P_j\cap S'$ is the maximal ideal of $S'$ and we have
 $$A_{\tau_j}\geq 1+\lceil\frac{\height (P_j\cap S')}{2}\rceil\geq 1+\lceil\frac{3}{2}\rceil\geq 3=\depth_S I.$$

If   $P_1+P_j\not =m\not =P_1+P_k$,$P_j+P_k\not =m$  then $$A_{\tau_{jk}}\geq \sdepth_{S_{\tau_{jk}}}J_{\tau_{jk}} +\sdepth_{S'}L_{\tau_{jk}}\geq 1+\sdepth_{S'}(P_j\cap P_k\cap S')\geq \depth_S I,$$
by Proposition \ref{1}.
If  $P_1+P_j =m\not =P_1+P_k$, $P_j+P_k\not =m$  then  $$\sdepth_{S'}L_{\tau_{jk}}\geq \depth_{S'}L_{\tau_{jk}}=1+\dim S/(P_1+P_k)\geq 1+q.$$ Thus $A_{\tau_{jk}}\geq 2+q\geq 3=\depth_SI$. If  $P_1+P_j =m=P_1+P_k$, $P_j+P_k\not =m$ then $L_{\tau_{jk}}$ is the maximal ideal of $S'$ and we get $A_{\tau_{jk}}\geq 1+\lceil\frac{3}{2}\rceil\geq 3=\depth_S I.$
If $I\cap S''\not=0$ then $A_0=\sdepth_S(I\cap S'')S\geq 1+n-r\geq \depth_SI$. Hence $\sdepth_S I\geq \depth_S I$.
\end{proof}

\vskip 1 cm

\section{Graph of a monomial squarefree  ideal}

 Let $I=\cap_{i=1}^s P_i$, $s\geq 3$  be the intersection of the minimal monomial prime ideals of $S/I$.    Assume that $\Sigma_{i=1}^s P_i=m$
and the bigsize of $I$ is two. We may suppose that $P_1=(x_1,\ldots,x_r)$ for some $r<n$ and set $$q=q(I)=\min\{\dim S/(P_i+P_j): j\not =i,P_i+P_j\not =m\}.$$ Thus $q\leq n-r.$ Set $S''=K[x_1,\ldots,x_r]$, $S'=K[x_{r+1},\ldots,x_n]$.

\begin{Definition}{\em
Let $\Gamma$ be the simple graph on $[s]$ given by $\{ij\}$ is an edge (we write $\{ij\}\in E(\Gamma)$) if and only if $P_i+P_j=m$.
 We call $\Gamma$ the {\em graph associated to} $I$.
 $\Gamma$ has the {\em triangle property}  if there exists $i\in [s]$ such that for all $j,k\in [s]$ with $\{ij\},\{ik\}\in E(\Gamma)$ it follows $\{jk\}\in E(\Gamma)$.
In fact the triangle property  says that it is possible to find  a "good" main prime in the terminology of  \cite[Example 4.3]{P2}, which we remind shortly next.}
\end{Definition}

\begin{Example}\label{vechi}{\em Let $n=10$, $P_1=(x_1,\ldots,x_7)$, $P_2=(x_3,\ldots,x_8)$,

\noindent $P_3=(x_1,\ldots,x_4,x_8,\ldots,x_{10})$, $P_4=(x_1,x_2,x_5,x_8,x_9,x_{10})$,
$P_5=(x_5,\ldots,x_{10})$ and $I=\cap_{i=1}^5 P_i$. Then $q(I)=2$, and $\depth_SI=4$. The graph associated to $I$  on [5] as above has edges $$E(\Gamma)=
\{\{13\},\{15\},\{35\},\{14\},\{23\},\{24\}\}$$
 and has the triangle property, but only $\{5\}$ is a "good" vertex, that is for all $j,k\in [4]$ with $\{j5\},\{k5\}\in E(\Gamma)$ it follows $\{jk\}\in E(\Gamma)$.
Below you have the picture of $\Gamma$.}
\end{Example}

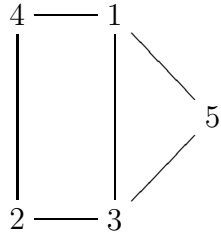
\begin{displaymath}
\xymatrix{4 \ar@{-}[r] \ar@{-}[dd] & 1\ar@{-}[rd] \ar@{-}[dd] & \\
            &&5\\
          2 \ar@{-}[r] & 3 \ar@{-}[ru]}
\end{displaymath}

\begin{Proposition}\label{2} If the bigsize of $I$ is two and $\Gamma=\Gamma(I)$ has  the triangle property then $\sdepth_S I\geq \depth_S I$.
\end{Proposition}
\begin{proof}
Renumbering $(P_i)$ we may suppose that $i=1$, that is for all $j,k\in [s]$ with $\{1j\},\{1k\}\in E(\Gamma)$ it follows $\{jk\}\in E(\Gamma)$ by the triangle property. We repeat somehow the proof of Proposition \ref{dep3}.
Applying \cite[Theorem 1.5]{P2} for ${\mathcal F}$ containing some $\tau_j=\{ j\}$,  and $\tau_{jk}=\{ j,k\}$  $1<j,k\leq s$,$j\not=k$
 we get $\sdepth_S I\geq \min\{A_0,\{A_{\tau_j}\}_{\tau_j\in {\mathcal F}}\},\{A_{\tau_{jk}}\}_{\tau_{jk}\in {\mathcal F}}\}$.
 Note that the bigsize of $J_{\tau}$ is $\leq 1$ (similarly $L_{\tau}$), $\tau\in {\mathcal F}$ and so 
  $\sdepth_{S_{\tau}}J_{\tau}\geq \depth_{S_{\tau}}J_{\tau}$ by \cite[Corollary 1.6]{P2}.
If $P_1+P_j\not =m$ then $$A_{\tau_j}\geq \sdepth_{S_{\tau_j}}J_{\tau_j} +\sdepth_{S'}(P_j\cap S')\geq 1+\dim S/(P_1+P_j)+\lceil\frac{\height (P_j\cap S')}{2}\rceil.$$
Thus $A_{\tau_j}\geq 2+q\geq \depth_S I$ by Lemmas \ref{g2}, \ref{b2}. If $P_1+P_j =m$ but there exists $e\not =j$ such that $P_e+P_j\not =m$, then $\sdepth_{S_{\tau_j}}J_{\tau_j}\geq \depth_{S_{\tau_j}}J_{\tau_j}=1+\depth_S S/(\cap_{u\not=j}(P_u+P_j))\geq 1+q$ and so again $A_{\tau_j}\geq 2+q\geq \depth_S I$. If $P_e+P_j=m$ for all $e\not =j$ then $\depth_S I=2$ by \cite[Lemma 1.2]{P2}
and clearly $A_{\tau_j}\geq \depth_S I$.

Now note that if $P_1+P_j\not =m\not =P_1+P_k$, $P_j+P_k\not =m$  then $$A_{\tau_{jk}}\geq \sdepth_{S_{\tau_{jk}}}J_{\tau_{jk}} +\sdepth_{S'}L_{\tau_{jk}}\geq 1+\sdepth_{S'}(P_j\cap P_k\cap S')\geq \depth_S I$$
by Proposition \ref{1}. If  $P_1+P_j =m$, $P_j+P_k\not =m$  then $P_1+P_k\not =m$ by the triangle property and $\sdepth_{S'}L_{\tau_{jk}}\geq \depth_{S'}L_{\tau_{jk}}=1+\dim S/(P_1+P_k)\geq 1+q$. Thus $A_{\tau_{jk}}\geq 2+q\geq \depth_SI$.
If $I\cap S''\not=0$ then as in the proof of Proposition \ref{dep3} $A_0\geq \depth_SI$. Hence $\sdepth_S I\geq \depth_S I$.
\end{proof}

\begin{Definition}\label{conca}
{\em The graph $\Gamma$ is a {\em join graph} if it is a join of two of its subgraphs, that is after a renumbering of the vertices  there exists $1\leq c<s$ such that $\{ij\}\in E(\Gamma)$ for all $1\leq i\leq c$, $c<j\leq s$. Thus in Lemmas \ref{g2}, \ref{b2} one may say that $\depth_S\ S/I =1$ if and only if the associated graph of $I$ is a join graph.
Let $\Gamma_1$, $\Gamma_2 $ be  graphs on $[r]$, respectively $\{r,r+1,\ldots,s\}$ for some integers $1<r\leq s-2$. Let $\Gamma$ be the graph on $[s]$ given by
$E(\Gamma)= E(\Gamma_1)\cup E(\Gamma_2)\cup\{\{ij\}: 1\leq i< r,r\leq j\leq s\}.$
We call $\Gamma$ the graph given by {\em concatenation of} $\Gamma_1$, $\Gamma_2$ {\em in the vertex} $\{r\}$.}
\end{Definition}

\begin{Lemma}\label{depthgr} Let $I=\cap_{i=1}^s P_i$ be the intersection of the minimal monomial prime ideals of $S/I$, $I_1=\cap_{i=1}^r P_i$,
$I_2=\cap_{i\geq r}^s P_i$ and  $\Gamma$,  $\Gamma_1$, $\Gamma_2 $  be the graphs associated to $I$, respectively $I_1$, $I_2$ as in the previous section. Suppose that $\Sigma_{i=1}^s P_i=m$,
$\Gamma$ is the concatenation of  $\Gamma_1$, $\Gamma_2 $ in $\{r\}$ and $bigsize(I)=2$. Then
$$\depth_SI=\min\{\depth_SI_1,\depth_SI_2\}.$$
\end{Lemma}
\begin{proof} We consider the following exact sequence
$$0\rightarrow S/I\rightarrow S/I_1\oplus S/I_2 \rightarrow S/I_1+I_2\rightarrow 0.$$
Since $P_i+P_j=m$ for all  $1\leq i<r$, $r<j\leq s$ we get $I_1+I_2=P_r$. But $\depth_S S/I,\depth_SS/I_i\leq \depth_SS/P_r$ for $i=1,2$ and by Depth Lemma \cite{Vi} we get
$$\depth_SS/I=\min\{\depth_SS/I_1,\depth_SS/I_2\}.$$
\end{proof}

\begin{Remark}\label{ap}  {\em Let  $I=\cap_{i=1}^3 P_i$ be the  intersection of the minimal monomial prime ideals of $S/I$. Suppose that $P_1+P_2\not =m\not =P_1+P_3$ and
$P_2+P_3=m$. Let  $I_1=P_1\cap P_2$,
$I_2=P_1\cap P_3$ and  $\Gamma$,  $\Gamma_1$, $\Gamma_2 $  be the graphs associated to $I$,  respectively $I_1$, $I_2$. We have $E(\Gamma_1)=E(\Gamma_2)=\emptyset$ and $E(\Gamma)=\{\{23\}\}$. Then $\Gamma$ is the concatenation of   $\Gamma_1$, $\Gamma_2 $ in $\{1\}$ and
$$\depth_S\ I=\min\{\depth_S\ I_1,\depth_S\ I_2\}=2+\min\{\dim\ S/(P_1+P_2),\dim\ S/(P_1+P_3)\}$$
by the above lemma. That is the formula found in \cite[Proposition 2.1]{AP}.}
\end{Remark}
\begin{Proposition}\label{descom} Let $I=\cap_{i=1}^s P_i$ be the intersection of the minimal monomial prime ideals of $S/I$,
and  $\Gamma$  be the graph associated to $I$. Suppose that $\Sigma_{i=1}^s P_i=m$, $bigsize(I)=2$ and $\depth_SI>3$. Then after renumbering $(P_i)$ there exists $1<r\leq s-2$ such that $\Gamma$ is the concatenation in $\{r\}$ of the graphs  $\Gamma_1$, $\Gamma_2$ associated to $I_1=\cap_{i=1}^r P_i$, respectively $I_2=\cap_{i\geq r}^s P_i$. Moreover, $\depth_S\ I_1,\depth_S\ I_2>3$.
\end{Proposition}
\begin{proof}
Since $bigsize(I)=2$ we may suppose that $P_{s-1}+P_s\not =m$, that is $\{s-1,s\}\not\in\Gamma$.  Consider the following exact sequence
$$ 0\rightarrow S/I \rightarrow S/P_1\cap \ldots \cap P_{s-1} \oplus S/P_1\cap \ldots \cap P_{s-2}\cap P_s\rightarrow S/P_1\cap \ldots \cap P_{s-2}\cap (P_s+P_{s-1})
\rightarrow 0.$$
As in the proof of Proposition \ref{2} we see that $P_i\not\subset P_s+P_{s-1}$ for $i<s-1$ because $bigsize(I)=2$. Then $\depth_SS/P_1\cap \ldots \cap P_{s-2}\cap (P_s+P_{s-1})=1$ using \cite[Corollary 1.6]{P2}. By Depth Lemma we get, let us say, $\depth_SS/P_1\cap \ldots \cap P_{s-1}=1$ since $\depth_SS/I>2$. After a renumbering of $(P_i)_{i<s-1}$ using Lemmas \ref{g2}, \ref{b2} we may suppose that there exists $1\leq c<s-1$ such that $P_i+P_j=m$ for all $1\leq i\leq c$, $c<j<s$.
Set $r=c+1$ and renumber $P_s$ by $P_r$ and $P_i$ by $P_{i+1}$ for $r\leq i< s$. Then $I_1=\cap_{i=1}^r P_i$ and $I_2=\cap_{i\geq r}^s P_i$ satisfy our proposition. The last statement follows by Lemma \ref{depthgr}.
\end{proof}

\begin{Corollary} \label{carac} Let $I=\cap_{i=1}^s P_i$ be the intersection of the minimal monomial prime ideals of $S/I$. Suppose that $\Sigma_{i=1}^s P_i=m$, $bigsize(I)=2$ and $q(I)>1$. Then $\depth_S\ I=2+q(I)>3$ if and only if the graph associated to $I$ is a concatenation of several graphs on two vertices having no edges.
\end{Corollary}
\begin{proof} The necessity follows  applying the above proposition by recurrence and the sufficiency follows applying Lemma \ref{depthgr} by recurrence.
\end{proof}

\begin{Corollary} \label{car} Let $I=\cap_{i=1}^s P_i$, $I'=\cap_{i=1}^s P'_i$ be the intersection of the minimal monomial prime ideals of $S/I$, respectively $S/I'$. Suppose that $\Sigma_{i=1}^s P_i=m=\Sigma_{i=1}^s P'_i$, $bigsize(I)=bigsize(I')=2$ and $q(I)=q(I')$.  If  the graphs associated to $I$, respectively $I'$ coincide,
then $\depth_S\ I=\depth_S\ I'$.
\end{Corollary}
\begin{proof} By the above corollary  $\depth_S\ I>3$ and  $\depth_S\ I'>3$ hold if and only if the graphs $\Gamma(I)$, $\Gamma(I')$    are  concatenations of several graphs on two vertices having no edges. Since  $\Gamma(I)=\Gamma(I')$ we get that $\depth_S\ I>3$ if and only if $\depth_S\ I'>3$.  But Lemmas \ref{g2}, \ref{b2} says that in this case $\depth_S\ I=2+q(I)=2+q(I')=\depth_S\ I'$. Note that  $\depth_S\ I=2$ holds if and only if $\Gamma(I)=\Gamma(I')$ is a join graph which happens if and only if  $\depth_S\ I'=2$. Then necessary $\depth_S\ I=\depth_S\ I'$ also in the  left case $\depth_S\ I=3$.
\end{proof}

\begin{Theorem} \label{char} The depth of a monomial squarefree ideal $I$ of $S$ such that  the sum of every three different of its minimal prime ideals is the maximal ideal of $S$, or more general a constant ideal of $S$. Then the depth of $I$ does not depend on the characteristic of $K$.
\end{Theorem}
\begin{proof} It is enough to suppose $\Sigma_{i=1}^s P_i=m$ and $bigsize(I)=2$, $size(I)=1$ by Corollary \ref{ind}. By Lemmas \ref{g2}, \ref{b2} (see also Remark \ref{conca}) $\depth_S\ I=2$ if and only if the graph $\Gamma(I)$ associated to $I$  is a join graph which is a combinatorial characterization and so does not depend on $p=$char $K$. By Corollary \ref{carac} $\depth_S\ I=2+q(I)>3$ if and only if $q(I)>1$ and $\Gamma(I)$ is a concatenation of several graphs on two vertices having no edges, the exact value of $\depth_S\ I$ being given by $q(I)$. Thus again $\depth_S\ I$ does not depend on $p$. Finally, $\depth_S \ I=3$  happens only when we are not in the above cases, that is,  it  does not depend on $p$.
\end{proof}

\begin{Remark} \label{f}{\em
The above theorem fails if just the sum of every four minimal prime ideals of $I$ is the maximal ideal of $S$. \cite[Examples 1.3 ii)]{HPV} says that the Stanley-Reisner ideal $I$ of the simplicial complex associated to the canonical triangulation of the real projective plane has $bigsize(I)=3$, $size(I)=2$ and $\depth_S\ I=4$ if char $K\not=2$, otherwise $\depth_S\ I=3$.}
\end{Remark}

\vskip 1 cm
\section{Stanley's Conjecture for monomial squarefree ideals of bigsize $2$}

 The case when $\depth_SI>1+bigsize(I)$ is unusual big and it is hard to check the Stanley's Conjecture in this case. Next we will construct such examples, where Lemma \ref{depthgr} and Proposition \ref{descom}  prove to be very useful.

 \begin{Example}\label{ex1} {\em  Let $\Gamma_1$ be the graph given on $\{1,2,5\}$ by $E(\Gamma_1)=\{\{15\}\}$ and
  $\Gamma_2$ be the graph given on $\{3,4,5\}$ by $E(\Gamma_1)=\{\{35\}\}$. Suppose that $I_1=P_1\cap P_2\cap P_5$ and $I_2=P_3\cap P_4\cap P_5$ are  irredundant intersections of monomial prime ideals of $S$ with $q(I_1)>1$, $q(I_2)>1$, $bigsize(I_1)=bigsize(I_2)=2$. Then $\depth_S\ I=2+\min\{\dim S/P_1+P_2,\dim S/P_2+P_5\}\geq 2+q(I_1)>3$ by \cite{AP} (see Remark \ref{ap}) if $q(I_1)>1$. Similarly, $ \Gamma_2$ is the graph associated to $I_2$ and $\depth_SI_2>3$. Let $\Gamma$ be the concatenation in $\{5\}$ of $\Gamma_1$ and $\Gamma_2$. If $I=I_1\cap I_2$ is an irredundant intersection of those $5$-prime ideals and $q(I)>1$, $bigsize(I)=2$ then $\Gamma$ is the graph associated to $I$  and $\depth_SI>3$ by Lemma \ref{depthgr}. This is the graph from the Example \ref{vechi}.}
  \end{Example}
  The above example is not bad because this case can be handled by our Proposition \ref{2}, that is there exists a "good" main prime $P_5$. Are there ideals $I$ for which there exists no "good" main prime, that is the graph associated to $I$ does not have the triangle property?  Next we will construct such a bad example. First we will see
  how should look its graph.

\begin{Example}\label{ex2} {\em Let $\Gamma_1$ be the graph constructed above on $[5]$ and $\Gamma_2$ be the graph given on $\{1,6\}$ with $E(\Gamma_2)=\emptyset$.
Let $\Gamma $ be the concatenation in $\{1\}$ of $\Gamma_1$ and $\Gamma_2$. Below you have the picture of $\Gamma$  and clearly it does not satisfy the triangle property. If we show that $\Gamma$ is the graph associated to a monomial squarefree ideal $I$ of $S$ with $bigsize(I)=2$ and $q(I)>1$ then we will have $\depth_SI>3$ by Lemma \ref{depthgr}.
This is done in the next example.}
\end{Example}
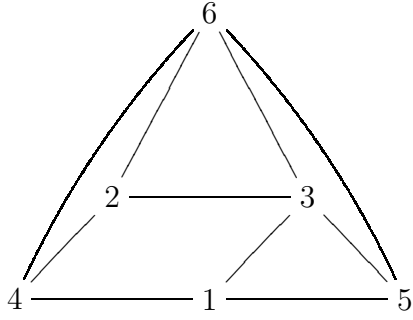
\begin{displaymath}
\xymatrix{ & & 6 \ar@{-}[ddl] \ar@/_/@{-}[dddll] \ar@{-}[ddr] \ar@/^/@{-}[dddrr] & &\\
            & & & & \\
            & 2 \ar@{-}[dl] \ar@{-}[rr] & & 3 \ar@{-}[dl] \ar@{-}[dr] &\\
            4 \ar@{-}[rr] & &1 \ar@{-}[rr] & & 5}
\end{displaymath}
\begin{Example}\label{ex3} {\em Let $n=12$, $P_1=(x_1,x_4,x_5,x_6,x_9,\ldots,x_{12})$, $P_2=(x_1,x_4,\ldots,x_{10})$, $P_3=(x_1,x_2,x_3,x_7,x_8,\ldots,x_{12})$,
$P_4=(x_1,x_2,x_3,x_6,x_7,x_8,x_{11},x_{12})$, $P_5=(x_1,\ldots,x_8)$, $P_6=(x_2,\ldots,x_6,x_9,\ldots,x_{12})$ and $I=\cap_{i=1}^6 P_i$. We have $P_1+P_4=P_1+P_5=
P_1+P_3=$
$$P_2+P_3=P_2+P_4=P_2+P_6=P_3+P_5=P_3+P_6=P_4+P_6=P_5+P_6=m$$ and
$P_1+P_2=m\setminus \{x_2,x_3\}$, $P_1+P_6=m\setminus \{x_7,x_8\}$, $P_2+P_5=m\setminus \{x_{11},x_{12}\}$, $P_3+P_4=m\setminus \{x_4,x_5\}$, $P_4+P_5=m\setminus \{x_9,x_{10}\}$. Clearly, $bigsize(I)=2=q(I)$ and the graph $\Gamma$ associated to $I$ is the graph constructed in Example \ref{ex2}. We have $\depth_S\ I=2+q(I)=4$ by Lemma \ref{depthgr}.
Let $S'=K[x_2,\ldots, x_{12}]$ and $P'_i=P_i\cap S'$, $I'=I\cap S'$. We have $I'=\cap_{i=1}^5 P'_i$ because $P'_1\subset P'_6$. The graph associated to $I'$ is in fact  $\Gamma_1$ from the above example and has the triangle property. Then by Proposition \ref{2} we get $\sdepth_{S'}\ I'\geq \depth_{S'}\ I'=2+q(I')=4$ because $q(I')=q(I)$.
Using the decomposition $I=I'\oplus x_1(I:x_1)$ as linear spaces we get
$$\sdepth_{S} I\geq \min\{\sdepth_{S'} I',\sdepth_{S} (I:x_1)\}\geq\min\{4,\sdepth_{S} P_6\}=$$
$$\min\{4,\lceil\frac{9}{2}\rceil+3\}=4=\depth_{S}\ I.$$
This gives us the idea to handle such bad examples in the next.}
\end{Example}

\begin{Proposition} \label{SC}Let $I=\cap_{i=1}^s P_i$ be the intersection of the minimal monomial prime ideals of $S/I$. Suppose that $\Sigma_{i=1}^s P_i=m$, $bigsize(I)=2$ and $\depth_S\ I>3$. Then $\sdepth_S\ I\geq \depth_S\ I$.
\end{Proposition}
\begin{proof}
Apply induction on $s$. The case $s\leq 4$ are given in \cite{PQ}, \cite{AP}, \cite{P2}. Assume that $s>4$ and let $\Gamma$ be the graph of $I$. By Proposition
\ref{descom} we may suppose after renumbering $(P_i)$ that there exists $1<r\leq s-2$ such that  $\Gamma$ is the concatenation in $\{r\}$ of the graphs  $\Gamma_1$,  $\Gamma_2$ associated to $I_1=\cap_{i=1}^r P_i$, respectively $I_2=\cap_{i=r}^s P_i$. Note that if $r=2$ or $s-r=2$  then $bigsize(I_i)$  could be not $2$ but this makes no troubles since we  need the bigsize  $2$ only to apply Proposition \ref{descom}. From Lemma \ref{depthgr} it follows that
$\depth_{S}\ I=\min\{\depth_{S}\ I_1,\depth_{S}\ I_2\}$ and so $\depth_{S}\ I_i>3$ for $i=1,2$. Note that $P_r\setminus P_j\subset P_i\setminus P_j=m\setminus P_j$ for all $1\leq i<r$, $r<j\leq s$. After renumbering variables we may suppose that $\{x_1,\ldots,x_e\}$ $1\leq e<n$ are all variables of $\cup_{j>r}^s (P_r\setminus P_j)$. As we noticed they are  contained in any $P_i$, $1\leq i<r$. Set $S'=K[x_{e+1},\ldots,x_n]$, $P'_i=P_i\cap S'$, $I'=I\cap S'$. Then  $P'_r\subset P'_j$ for all $r<j\leq s$ and we get $I'=\cap_{i=1}^r P'_i$. Moreover, since $\{x_1,\ldots,x_e\}$ is contained in any $P_i$, $1\leq i\leq r$ we see that the "relations" between these prime ideals preserve after intersection with $S'$ and the graph $\Gamma'$ of $I'$ is in fact  $\Gamma_1$. Moreover, $q(I')=q(I_1)$ and $bigsize(I')=bigsize(I_1)$. Then $\depth_{S'}\ I'=\depth_{S}\ I_1$ by Corollary \ref{car}, the case $r=2$ being trivial. Using induction hypothesis on $s$ we get
$\sdepth_{S'}\ I'\geq \depth_{S'}\ I'$. We have the decomposition $I=I'\oplus ((x_1,\ldots,x_e)\cap I)$  as linear spaces and it follows
$$\sdepth_{S}\ I\geq \min\{\sdepth_{S'}\ I', \sdepth_S\ ((x_1,\ldots,x_e)\cap I) \}.$$
But $J=(x_1,\ldots,x_e)\cap I=\cap_{i>r}^s P'_i\cap (x_1,\ldots,x_e)$ because $(x_1,\ldots,x_e)\subset P_i$ for $1\leq i\leq r$ and the decomposition is irredundant
since if $(x_1,\ldots,x_e)\subset P_j$ then $P_r\subset P_j$ which is false.  Note that $q(J)=q(I_2)$ and the graph associated to $J$ coincides with $\Gamma_2$.
Again  by Corollary \ref{car} $\depth_{S}\ J=\depth_{S}\ I_2>3$, the case $r=2$ being trivial. Using induction hypothesis on $s$ we get
$\sdepth_S\ J\geq \depth_S\ J$ and so
$$\sdepth_S\ I\geq \min\{\sdepth_{S'}\ I', \sdepth_S\ J \}\geq \min\{\depth_{S}\ I_1,\depth_{S}\ I_2\}=\depth_{S}\ I.$$
\end{proof}

\begin{Theorem} \label{m} Let $I$ be a  monomial squarefree ideal  of $S$ such that  the sum of every three different of its minimal prime ideals is the maximal ideal of $S$, or more general a constant ideal $J$ of $S$. Then $\sdepth_S\ I\geq \depth_S\ I$.
\end{Theorem}
The proof follows from Theorem \ref{dep3} and the above proposition, the reduction to the case $J=m$ being given by \cite[Lemma 3.6]{HVZ}.

\end{document}